\documentclass[12pt]{article}

\usepackage{amsmath}
\usepackage{amssymb}
\usepackage{amsthm}
\usepackage{MnSymbol}
\usepackage{lineno}

\newtheorem{thm}{Theorem}[section]
\def \btm {\begin{thm}}
\def \etm {\end{thm}}
\newtheorem{prp}[thm]{Proposition}
\def \bpn {\begin{prp}}
\def \epn {\end{prp}}
\newtheorem{lem}[thm]{Lemma}
\def \blm {\begin{lem}}
\def \elm {\end{lem}}
\newtheorem{obs}[thm]{Observation}
\def \bob {\begin{obs}}
\def \eob {\end{obs}}
\newtheorem{rmk}[thm]{Remark}
\def \brm {\begin{rmk}}
\def \erm {\end{rmk}}
\newtheorem{prm}[thm]{Problem}
\def \bpm {\begin{prm}}
\def \epm {\end{prm}}
\newtheorem{dfn}[thm]{Definition}
\def \bdf {\begin{dfn}}
\def \edf {\end{dfn}}
\newtheorem{exa}[thm]{Example}
\def \bex {\begin{exa}}
\def \eex {\end{exa}}
\newtheorem{exas}[thm]{Examples}
\def \bexs {\begin{exas}}
\def \eexs {\end{exas}}
\def \bpf {\begin{proof}}
\def \epf {\end{proof}}

\def \msk {\medskip}
\def \ssk {\smallskip}
\def \nin {\noindent}
\def \smin {\setminus}

\def \cF {\mathcal{F}}
\def \vp {\varphi}

\def \arr {anti-Ramsey}
\def \Arr {Anti-Ramsey}

\def \AR {\mathrm{ar}}
\def \arra {\mathrm{AR}}

\def \soar {\mathrm{Sod}}

\def \lar {\mathrm{Lr}}

\def \arb {\mathrm{Arb}}
\def \pla {\mathrm{Pla}}

\def \ex {\mathrm{ex}}

\def \erd {Erd\H os}
\def \kov {K\H ov\'ari}
\def \sos {S\'os}
\def \tur {Tur\'an}

\def \chiF {\chi_{_\cF}}
\def \chiFk {\chi_{_{\cF_k}}}
\def \chim {\chi^-}

\begin{document}

\title{Monochromatic graph decompositions \\
  inspired by \arr\ colorings}
\author{Yair Caro\,\thanks{~Department of Mathematics, University of Haifa-Oranim, Tivon 36006, Israel}
 \and Zsolt Tuza\,\thanks{~HUN-REN Alfr\'ed R\'enyi Institute of Mathematics, Budapest, Hungary} $^,$\thanks{~Department of Computer Science and Systems Technology, University
of Pannonia, Veszpr\'em, Hungary} $^,$\thanks{~Research supported in part by the National
Research, Development and Innovation Office, NKFIH Grant FK 132060}}
\date{~~}
\maketitle

\vspace{-7ex} 

\begin{abstract}
We consider coloring problems inspired by the theory of \arr\
 / rainbow colorings that we generalize to a far extent.  

Let $\cF$ be a hereditary family of graphs; i.e., if
 $H\in \cF$ and $H'\subset H$ then also $H'\subset \cF$.
For a graph $G$ and any integer
 $n \geq |G|$, let $f(n,G|\cF)$ denote the smallest number
 $k$ of colors such that any edge coloring of $K_n$ with
 at least $k$ colors forces a copy of $G$ in which each color
 class induces a member of $\cF$.  
 
The case $\cF = \{K_2\}$ is the notorious \arr\ /
 rainbow coloring problem introduced by \erd, Simonovits
 and \sos\ in 1973.  

Using the $\cF$-deck of $G$,
 $D(G|\cF) = \{  H : H = G - D, \, D \in \cF\}$,
  we define
  $\chiF(G) = \min \{ \chi(H) :  H \in D(G|\cF) \}$.

The main theorem we prove is: 
Suppose $\cF$ is a hereditary family of graphs, and
  let $G$ be a graph not a member of $\cF$.
 
(1)\quad If $\chiF(G) \geq 3$, then
 $f(n, G |\cF) = (1+o(1)) \, \ex(n, K_{\chiF(G)})$.
 
(2)\quad Otherwise $f(n, G |\cF)  =  o(n^2)$.
 
Among the families covered by this theorem are:
 matchings, acyclic graphs,
 planar and outerplanar graphs, $d$-degenerate graphs,
 graphs with chromatic number at most $k$, graphs with
 bounded maximum degree, and many more.

We supply many concrete examples to demonstrate the wide range
 of applications of the main theorem; the next result is a
  representative of these examples.   

For $p \geq 5$ and $\cF = \{ tK_2 : t \geq 1 \}$, we have
   $f(n,K_p |\cF)  = (1+o(1)) \, \ex(n, K_{\lceil p/2 \rceil})$;
 this means a properly colored copy of $K_p$.
In other words, a certain number of colors forces nearly twice as large
 properly edge-colored complete subgraphs as rainbow ones.
\end{abstract}

\section{Introduction}
\label{s:intro}

The main purpose of this paper is to offer a wide-range
 generalization of the rainbow anti-Ramsey theory initiated
  by \erd, Simonovits and \sos, 50 years ago \cite{ESS-75}.

\msk

An edge-colored graph $G$ is called rainbow if all its
 edges are colored with distinct colors. 
The anti-Ramsey number denoted by $\AR(n,G)$ for $n\geq |G|$
 is defined as the smallest $k$ such that, in any coloring 
of the edges of $K_n$ using at least $k$ colors, there is a rainbow 
copy of $G$.
We note that sometimes $\AR(n,G)$ is defined as the maximum $k$ 
for which rainbow $G$ can be avoided. 

The most general result in this direction was given already in \cite{ESS-75}.
Let $\chim(G) = \min \{\chi(G - e): e  \in  E(G)\}$.
Then $\AR(n,G) = (1+o(1)) \, \ex(n, K_{\chim(G)})$, 
where in general $\ex(n, H)$ is the celebrated \tur\ number for a graph $H$, 
and $\chi(H)$ is the chromatic number of $H$.  

Since then, hundreds of papers were written on this major
 coloring problem, and research in this direction is still
 active; see, e.g., the recent publications
  \cite{ADK-24+,A-etal-24,AC-24,B-etal-24+,H-24t,JLZ-19,JG-24,LSS-21,LTY-24+,LS-23+,XLYX-22}.  
For further information concerning \arr\ problems and alike,
  we refer to the detailed study of small graphs in \cite{BGR-15}
 and to  the dynamic survey \cite{dyn-surv}.

To present our generalization we need a couple of definitions,
 parallel to those in anti-Ramsey theory.

A family $\cF$ of graphs is hereditary if it is closed under
 inclusion of subgraphs; namely, if $G \in \cF$ and $H \subset G$
 then $H  \in \cF$.
We refer to \cite{hered-surv} for information concerning
 properties of hereditary families of graphs. 

For a given graph $G$ and hereditary\footnote{In some
 of the definitions it suffices to assume that $G$ contains
 at least one subgraph $D\in\cF$. In some others in the
 sequel we need that the edge set of $G$ is decomposable
 into edge-disjoint subgraphs that are members of $\cF$.
 A sufficient condition for both is $K_2\in\cF$.}
 family $\cF$ we define the $\cF$-deck of $G$ as 
  $$
    D(G|\cF) = \{H:  H = G - D, \ D  \in \cF \} \,,
  $$
 and the reduced chromatic number of $G$ 
 with respect to $\cF$ as
  $$
    \chiF(G) = \min \{\chi(H): H  \in D(G|F) \} \,.
  $$
The main problem we consider is:  determine or supply good
 approximation to $f(n, G | \cF)$, the smallest integer $k$
 such that, in every coloring of the edges of $K_n$ using at
 least $k$ colors, there is a copy of $G$ in which all edges
 of the same color induce a member of $\cF$. 

We assume here and in the sequel throughout $n \geq |G|$,
 to avoid the situation that a smaller $K_n$ trivially
 does not contain any copy of $G$.
Clearly the case of $F = \{K_2\}$ being a single edge
 constitutes the original anti-Ramsey theory.   

Our main result is as follows.

\btm
Let\/ $G$ be a graph and\/ $\cF$ a hereditary family of graphs.
 \begin{itemize}
   \item[$(i)$] If\/ $\chiF(G) \geq 3$, then\/
    $f(n, G|\cF) = (1+o(1)) \, \ex(n, K_{\chiF(G)})$.
   \item[$(ii)$] If\/ $\chiF(G) = 2$, then\/
    $f(n, G|\cF) = o(n^2)$.
   \item[$(iii)$] If\/ $\chiF(G) = 1$, then\/ $f(n, G|F) = 1$.
 \end{itemize}
\etm

A highly interesting particular case is where $\cF$ is the
 family of matchings, $\cF = \{tK_2 : t \geq 1\}$, and
 $G = K_p$ is the complete graph of any given order $p\geq 5$.
Here,
  $$
    f(n, K_p |\cF) = (1+o(1)) \, \ex(n,K_{\lceil k/2 \rceil})
  $$
 holds, meaning a properly edge-colored copy of $K_p$. 
In other words, a certain number of colors forces a nearly
 twice as large properly edge-colored complete subgraph
 as a rainbow one. 

Due to the high level of generality of our new model,
 the paradigm behind the above-mentioned main theorem 
 in anti-Ramsey theory can be extended to the
 much wider range of hereditary families of graphs.
However, there is a major distinction between the two approaches.
In anti-Ramsey theory the proof of the upper bound requires only a 
simple application of the celebrated \erd--Stone--Simonovits
 theorem \cite{ES-66,ES-46} to find a large rainbow
 complete multipartite subgraph, and then
 add just one further edge with a distinct color.  
This is far from being the case in our more general setting.
Indeed, as $\chiF(G)$ can be much smaller than $\chi(G)$,
 we need to use many edges from the non-rainbow part
 to get the desired $\cF$-colored copy of $G$, preserving
 the required property that every edge-color-class induces
 a member of $\cF$.  
To deal with this more general setting we develop two lemmas,
 that allow us to overcome the obstacle explained above.
These lemmas will be presented in Section \ref{s:transv}.

\msk

The paper is organized as follows: 

In Section \ref{s:LB} we prove a general lower bound;
 namely, if $G$ is not a member of $\cF$, then
  $f(n, G|\cF) \geq \ex(n,K_{\chiF(G)}) +2$.

In Section \ref{s:transv} we develop an ``Independent
 Transversal Lemma'' and a ``Rainbow Cut Lemma'', which are
 the main tools---of potential further ap\-plications---needed
  for the proof of the main theorem. 

In Section \ref{s:UB} we prove the asymptotic upper bound, by
 combining the tools developed in Section \ref{s:transv} with
 the \erd--Stone--Simonovits theorem if $\chiF(G) \geq 3$, or
 with the \kov--\sos--\tur\ theorem \cite{KST-54}
  in case of $\chiF(G) = 2$. 

In Section \ref{s:applic} we provide concrete examples of the main theorem when $\cF$ is 
the family of: outerplanar graphs, planar graphs, acyclic graphs, 
matchings, $d$-degenerate graphs, $k$-colorable graphs, graphs of maximum degree at most $d$, 
while $G$ is 
either any graph or in some examples the complete graph $G = K_p$. 
Also, we give an asymptotic solution of the sub-$k$-coloring
 version of the \arr\ problem, studied recently in \cite{H-24t},
 for all graphs that cannot be made bipartite by the removal
 of at most $k$ edges.

In Section \ref{s:concl} we shall offer some further problems
 for future research, and analyze the tightness of the estimate
 given on independent transversals in Section \ref{s:transv}. 

\msk

We mention here that in a companion paper \cite{ar-M} we consider 
further problems of this kind where the inspiration comes from a discussion 
with Riste \v Skrekovski concerning odd-coloring \cite{CPS-22} and 
conflict-free colorings \cite{CPS-23}. 
However, the main distinction is that in these ``odd-coloring variants'' 
or ``conflict-free variants'' the ``coloring family'' $\cF$ is non-hereditary, 
and the required tools are different.

\paragraph{Notation and Definitions.}

We follow standard notation as in West \cite{W-intro}.
In particular, we write $|G|$ for
 the order of $G$ (the number of its vertices); and for a
 subgraph $H$ of $G$, we denote by $G-H$ the graph obtained
 from $G$ by deleting the edge set $E(H)$ of $H$.
We shall add further necessary notation and definitions at
 later points.

\section{Reduced chromatic number, general lower bound, and non-critical subgraphs}
   \label{s:LB}

A decomposition of a graph
    $G=(V,E)$ is a collection of subgraphs
    $G_1=(V_1,E_1),\dots,G_m=(V_m,E_m)$ where
     $E_1\cup\cdots\cup E_m=E$ and the $G_i$ are
    mutually edge-disjoint.
For a given family $\cF$ of graphs,
    an $\cF$-decomposition requires
 $G_i\in\cF$ for all $1\leq i\leq m$.

Note that every graph $G$ admits an $\cF$-decomposition
 whenever $\cF$ is hereditary.
In fact, for decomposability it suffices to assume that
 the single edge $K_2$ belongs to $\cF$.
Given $\cF$ and $G$, let us denote by $F(G)$ the minimum number
 of members of $\cF$ that form an $\cF$-decomposition of $G$.
We next consider hereditary families in connection with
 edge decompositions and variants of \arr\ functions.

\bex   \label{exs:F}
Some interesting natural graph classes are as follows.
 \begin{enumerate}
  \item $\cF=\{K_2\}$, a single edge:\/
   $F(G)=|E(G)|$,\/ $f(n,G |\cF) = \AR(n,G)$.
  \item $\cF=\{tK_2 \mid t\geq 1\}$,\/ matchings of all sizes:\/
   $F(G)=\chi'(G)$ (chromatic index),\/
    $f(n,G |\cF) = \lar(n,G)$; notation\/ $\lar$ abbreviates
    ``\,local rainbow'' coloring.

  \item $\cF =$ graphs with at most\/ $k$ edges \cite{H-24t}:\/
    $F(G) = \lceil\, |E(G)|/k \rceil$,\/ $f(n,G |\cF)= \arra_k(n,G)$.
  \item $\cF =$ graphs of maximum degree\/ $t$:\/  $f(n,G |\cF)= \lar(n,t,G)$.
  \item $\cF$ = planar graphs:\/ $F(G) =$ thickness of\/ $G$,\/
   $f(n,G|\cF) = \pla(n,G) =$ smallest\/ $k$ such that
    every edge coloring of\/ $K_n$ using at least\/ $k$ colors
    forces a copy of\/ $G$ where each color class in\/ $G$
    induces a planar graph. 
  \item $\cF =$ acyclic graphs:\/ $F(G) =$ arboricity,\/
   $f(n,G|\cF) = \arb(n,G) =$ smallest\/ $k$ such that
    every edge coloring
    of\/ $K_n$ using at least\/ $k$ colors forces a copy of\/ $G$
     where all edge color classes in\/ $G$ are acyclic\}.  
  \item $\cF =$ odd graphs:\/
    $f(n,G|\cF) = \soar(n,G)$; notation\/ $\soar$ abbreviates
    ``strong odd'' coloring.
 \end{enumerate}
\eex

\brm
There are some global \emph{versus} local analogies
 in the conditions\/ $1\leftrightarrow 2$ and\/
  $3\leftrightarrow 4$ of the above list:\/
 $\AR(n,G)$ completely excludes repetition of colors
  in a suitable copy of\/ $G$, whereas
  $\lar(n,G)$ excludes repetition at any vertex in that copy.
Similarly,\/ $\arra_k(G,K_n)$ bounds the number of occurrences
 of each color in the entire\/ $G$, whereas\/ $\lar(n,t,G)$
 puts an upper bound at each vertex of\/ $G$.
\erm

We are now ready to present our main general lowerbound.

\btm   \label{t:LB}
Let\/ $\cF$ be a hereditary family, and\/ $G\notin\cF$ any graph.
 Then\/ $f(n,G, |F) \geq \ex( n,D(G|\cF)) +2$.
\etm

\bpf
Let $R$ be a rainbow extremal graph for $D(G | \cF)$
 with $\ex(n,D(G | \cF))$ edges, and assign a new color $c$
 to all the other edges of $K_n$.
 
Suppose for a contradiction that there is an $\cF$-colored copy of $G$. 
Since no member of $D(G|\cF)$ is in $R$, $G$ itself is not a subgraph of $R$.
But then the $\cF$-colored copy of $G$ must use at least one $c$-colored edge.
However, as $G$ is not a member of $\cF$, $G$ cannot use only $c$-colored edges and must use also edges from $R$.
 
Delete the $c$-colored edges (forming a member of $\cF$)  from $G$. 
The resulting graph $H$ is a member of $D(G |\cF)$, but
 by the construction of $R$, it contains no member of
 $D(G|\cF)$, a contradiction proving the assertion.
\epf

Next, we discuss the important role of ``critical subgraphs''.
Recall that in the \arr\ setting
 $\chim(G) := \min \{ \chi(G-e) : e \in E(G) \}$ is 
 crucial in determining $\AR(n,G)$ asymptotically; and
 if $G$ is non-bipartite and $\chim(G) = \chi(G)$, i.e.\ $G$
  contains no critical edge with respect to the chromatic
 number, then $\AR(n,G) = (1+o(1)) \, \ex(n, K_{X(G)})$ holds.

In analogy with this, we call a subgraph $D \in \cF$
 of $G$ critical with respect to $\cF$ if
 $\chi(G - D) < \chi(G)$, hence forcing $\chiF(G) < \chi(G)$.
Further, we say that $G$ is stable with respect to $\cF$,
 if $\chiF(G) = \chi(G)$; that is, no subgraph $D\subset G$
 with $D\in \cF$ is critical with respect to $\cF$.

The relevance of stability concerning $f(n,G|\cF)$ is
 demonstrated by the following tight asymptotic result.

\btm   \label{t:stable}
Let\/ $\cF$ be a hereditary family of graphs, and\/ $G$ a
 stable graph with respect to\/ $\cF$.
If\/ $\chiF(G)\geq 3$, then\/
 $f(n, G|\cF) = (1+o(1)) \, \ex(n,K_{\chi(G)})$.
\etm

\bpf
Since $G$ is stable, namely $\chiF(G) = \chi(G)$,
 the lower bound $f(n, G|\cF) > \ex(n,K_{\chi(G)})$
  folows by Theorem \ref{t:LB}.

For an upper bound with the same asymptotics
 recall from \cite{ES-66,ES-46}
 that the \tur\ number of the complete $\chi(G)$-partite graph
  $K=K_{|G|,\dots,|G|}$ with all partite sets having $|G|$
  vertices is equal to $(1+o(1)) \, \ex(n,K_{\chi(G)})$.
Hence, if $\ex(n,K)+1$ or more colors are used, we can select
 one edge from each of the first $\ex(n,K)+1$ color
 classes---which is still $(1+o(1)) \,
  \ex(n,K_{\chi(G)})$---and in this way find a rainbow copy
 of $K$, which obviously contains a rainbow subgraph $G$.
\epf

In Section \ref{s:applic} we present several consequences of
 this theorem: concerning planar graphs, outerplanar graphs,
 and graphs of large girth.

\brm
The condition ``\,hereditary'' in Theorems \ref{t:LB} and
 \ref{t:stable} can be weakened to\/ $K_2\in\cF$.
This fact makes the results relevant for\/ $\soar(n,G)$,
 among many further cases of\/ $\cF$,
 where the family of odd graphs is not hereditary but still
 contains all matchings.
The approach is applicable to obtain both lower bounds and
 asymptotics, in case of stability with respect to\/ $\cF$.
\erm

\section{Independent transversals and rainbow cuts}
\label{s:transv}

In this section we develop the main tools required for
 the proofs of upper bounds on $f(n, G | \cF)$.

The next result for directed graphs can be thought of as
 complementary to independent transversals in undirected
  graphs---a much studied problem, see the recent reference \cite{GS-22,HW-23+,WW-22}.
We should note that the assumptions below allow arbitrarily
 large degrees in the undirected underlying graphs even
 locally, in the union of just two vertex classes also,
 which is not the case in the related results
 known for undirected graphs.

\btm {\bf\itshape (Independent Transversal Lemma)}   \label{l:ITL}
 \
Let\/ $H$ be a directed graph with vertex equipartition\/
 $V=V_1\cup\cdots\cup V_m$ where\/ $|V_i|=s$ for\/ $i=1,\dots,m$,
 and let the maximum out-degree of\/ $H$ be\/ $\Delta\!^+$.
\begin{itemize}
 \item[$(i)$] If\/ $s>m\Delta\!^+$, then there exists an
  independent set $T$ that meets each $V_i$.
 \item[$(ii)$] If\/ $s\geq (2r+m)\Delta\!^+ + r$,
  then there exists an independent set $T$ that shares
  at least $r$ vertices with each $V_i$.
\end{itemize}
\etm

\bpf

$(i)$\quad
The number of sets meeting each $V_i$ in exactly one vertex
 is equal to $s^m$.
Any arc between some $V_i$ and $V_j$ is contained in
 precisely $s^{m-2}$ of those sets.
The digraph contains at most $ms\Delta\!^+$ arcs.
Thus, if $s>m\Delta\!^+$, then a required set $T$ exists that
 does not contain any arcs of $H$.

$(ii)$\quad
Repeatedly applying $(i)$ $k$ times, we obtain that
 $s\geq k+m\Delta\!^+$
 implies the existence of $k$ mutually disjoint independent
 sets $T_1,\dots,T_k$, each of them meeting each $V_i$.
Now we construct an auxiliary graph with vertex set
 $x_1,\dots,x_k$, where $x_{j_1}x_{j_2}$ is an edge if and
 only if there is an arc between $T_{j_1}$ and $T_{j_2}$ in $H$.
The average degree is at most $2\Delta\!^+$, hence \tur's
 theorem yields that the independence number is at least
 $\frac{k}{2\Delta\!^+ + 1}$.
The union of $r$ sets $T_j$ corresponding to an independent set
 of size $r$ is a set $T$ with the required properties.
Consequently, a sufficient condition is
 $k\geq (2\Delta\!^+ + 1)r$, what means
 $s\geq (2r+m)\Delta\!^+ + r$.
\epf

Tightness of part $(i)$ above will be discussed in Section \ref{s:concl}.

\ssk

For asymptotic estimates on $f(n,G|\cF)$ the following result will be crucial.

\btm {\bf\itshape (Rainbow Cut Lemma)}   \label{l:RCL}
 \
For any positive integers\/ $m,p$ there exists\/ $q=q(m,p)$
 with the following property.
If\/ $\psi$ is an edge coloring of\/ $K_{qm}$,
 and\/ $K=K_{q,\dots,q}$ is a rainbow spanning subgraph of\/
 $K_{qm}$ under\/ $\psi$, then $K$ contains a (rainbow)
 subgraph\/ $K'=K_{p,\dots,p}$ of order\/ $pm$
 such that the colors of edges
 inside the vertex classes of\/ $K'$ in\/ $K_{qm}$
 do not occur in\/ $K$.
\etm

\bpf
Write $q$ in the form $q=sp$, where $s$ will be chosen later.
We partition each vertex class $V_i$ of $K$ into $s$ sets
 $V_{i,j}$ of size $p$ each.
Note that inside each $V_{i,j}$ at most $\binom{p}{2}$ colors
 can occur, and each such color appears in $K$ at most once.
We construct a digraph with $sm$ vertices $x_{i,j}$ that
 represent the sets $V_{i,j}$.
There is an arc from $x_{i_1,j_1}$ to $x_{i_2,j_2}$ if
 $V_{i_1,j_1}$ contains an edge whose color appears on some
 edge incident with $V_{i_2,j_2}$ in $K$.

This digraph has maximum out-degree
 $\Delta\!^+ \leq 2\binom{p}{2}$.
By Theorem \ref{l:ITL} $(i)$, in case of $s>m\Delta\!^+$
 there exist $m$ sets $X_1=V_{i_1,j_1},\dots,X_m=V_{i_m,j_m}$
 of size $p$ inside  $V_1,\dots,V_m$, respectively,
 such that the colors inside those $X_i$ do not occur
 between any two of those sets.
Consequently $X_1,\dots,X_m$ can be taken as the vertex
 classes of $K'$.

It follows that $s=mp^2$ and $q=mp^3$ are suitable choices.
\epf

\section{Asymptotic upper bounds}   \label{s:UB}

Our main goal is to prove that the lower bound in
 Theorem \ref{t:LB} is asymptotically tight whenever $\chiF(G)$
 is at least 3.
Before that, we give a subquadratic upper bound for
 the other cases.

\vbox{
\btm   \label{t:subex}
Let\/ $\cF$ be a hereditary family and\/
 $G$ a graph with\/ $\chiF(G)\leq 2$.
\begin{itemize}
 \item[$(i)$] If\/ $\chiF(G)=1$, then\/ $f(n,G |\cF)=1$.
 \item[$(ii)$] If\/ $\chiF(G)=2$, then\/ $f(n,G |\cF)=o(n^2)$.
\end{itemize}
\etm
}

\bpf
The condition $\chiF(G)=1$ is equivalent to $G\in\cF$ and,
 since $\cF$ is hereditary, also every subgraph of $G$ is
 a member of $\cF$.
Hence, no matter how many colors are used in $K_n$,
 in every copy of $G$ each color class induces a graph
 belonging to $\cF$.
This proves $(i)$.

To prove $(ii)$, suppose $H \in D(G |\cF)$ and $\chi(H) = 2$.  
Let $A \cup B  = V(H)$ be a bipartition of $H$.
This $H$ is obtained by the removal of $D\in\cF$ with
 no edge between its two parts $D_1 = G[A]$ and $D_2 = G[B]$.

We now choose a large $t$, and consider any edge coloring
 $\psi$ of $K_n$ with more than $\ex(n,K_{t,t})$ colors.
Picking one edge from each color class, we obtain a
 $K=K_{t ,t}$ rainbow-colored under $\psi$.
According to the Rainbow Cut Lemma \ref{l:RCL}, if $t$ is
 large enough with respect to $|G|$, we can guarantee $D_1$
 on one side and $D_2$ on the other side of $K$, such that
 no color of $D_1\cup D_2$ occurs in $K$.
But then we can supplement $D_1\cup D_2$ with edges from $K$
 and obtain an $\cF$-colored copy of $G$.
Consequently, if $\psi$ uses $\ex(n,K_{t,t})+1$ colors
 at least, then a required $G$ occurs.
Due to the \kov--\sos--\tur\ theorem \cite{KST-54} we have
 $\ex(n,K_{t ,t}) \leq O(n ^{2-1/t }) = o(n^2)$,
 implying the asserted upper bound.
\epf

\btm   \label{t:ex}
Let\/ $\cF$ be a hereditary family and\/
 $G$ a graph with\/ $\chiF(G)\geq 3$.
Then
   $$
     f(n,G |\cF) = (1+o(1)) \, \ex(n,K_{\chiF(G)}) \,.
   $$
\etm

\bpf
The assertion was already proved in Theorem \ref{t:stable}
 for all $G$ that are stable with respect to $\cF$.
So, here we may assume that $\chiF(G)<\chi(G)$ holds.
We will use the notation $m:=\chiF(G)$ and $p:=|G|$.
By assumption we have $m\geq 3$.

We know from Theorem \ref{t:LB} that $\ex(n,K_m)$ is
 a lower bound on $f(n,G |\cF)$.
For a matching asymptotic upper bound we select
 one $D\in\cF$ such that $\chi(G-D)=m$ and, under
 this condition, $|E(D)|$ is smallest.
Let $V_1,\dots,V_m$ be the color classes of a (arbitrarily
 chosen) proper $m$-coloring of $G-D$, and denote by $D_i$
 the subgraph of $D$ induced by $V_i$.
Here the minimality condition on $|E(D)|$ implies
 $D=D_1\cup\cdots\cup D_m$ because any subgraph of $D\in\cF$
 is also in $\cF$, and deleting just the edges of all $D_i$
 the chromatic number of $G$ already goes down to $m$.

With reference to Theorem \ref{l:RCL}, we let $q$ denote the
 integer $q(m,p)$ that guarantees the conclusion of the
 Rainbow Cut Lemma, informally stating that in any edge
 coloring of $K_n$, a large complete multipartite rainbow
 subgraph contains a still fairly large one that avoids all
 the colors occurring inside its vertex classes.

The \erd--Stone--Simonovits theorem \cite{ES-46,ES-66} states
 for any fixed $m\geq 3$ and any $q$ that, as the number $n$
  of vertices gets large, $(1+o(1)) \, \ex(n,K_m)$ edges are
 enough to guarantee the presence of not only $K_m$ but also
 a complete $m$-partite subgraph $K^*=K_{q,\dots,q}$.

Let $\psi$ be any edge coloring of $K_n$ with at least
 $\ex(n,K^*)+1$ colors.
We select one edge from each color class and find
 a rainbow $K^*$.
Theorem \ref{l:RCL} with the current choice of $q$ yields a
 complete $m$-partite rainbow $K'=K_{p,\dots,p}\subset K^*$,
 say with vertex set $X_1\cup\cdots\cup X_m$, such that $K'$
 does not contain any colors induced by any $X_i$.

Since $|X_i|=p$, there exists an embedding $\eta$ of $G$ into
 $K_n$ such that $\eta(V_i)\subset X_i$ for all $i=1,\dots,m$.
Let us denote $Y_i:=\eta(V_i)$.
The set $Y_1\cup\cdots\cup Y_m$ restricts $K'$ to an
   $m$-partite rainbow graph $K$ on $p$ vertices, and
  also specifies a complete subgraph $K_p\subset K_n$.
Here we have $\eta(G-D)\subset K$, and $\eta(D)\subset K_p-K$.
The color sets of these two parts are disjoint, hence each
 color occurring in $\eta(D)$ belongs to a color class of
 $\eta(G)$ which is entirely contained in $\eta(D)$.
Consequently those color classes are isomorphic to subgraphs
 of $D\in \cF$, and therefore members of $\cF$,
  because $\cF$ is hereditary;
 and the color classes in $\eta(G-D)$ are single edges.
In this way a copy of $G$ is found where
 each color class induces a graph from $\cF$.
This completes the proof.
\epf

\brm
Our approach offers an alternative proof of the
 \erd--Si\-mo\-no\-vits--\sos\ theorem, as follows.
If\/ $\chi(G-e)=\chi(G)$ holds for all edges\/ $e\in E(G)$,
 then Theorem \ref{t:stable} on stable graphs does the job.
Otherwise assume that\/ $e$ satisfies\/ $\chi(G-e)=\chi(G)-1$.
If\/ $\psi$ is an edge coloring with more colors than the
 claimed upper bound, we select one edge from each color
 class and find a rainbow complete\/ $(\chi(G)-1)$-partite
 $K=K_{|G|,\dots,|G|}$.
Add any one edge\/ $e$ from a vertex class of\/ $K$, and
 remove the only one possible vertex\/ $v$ of\/ $K$ that is
 connected to\/ $e$ with a color\/ $\psi(e)$.
Still there is more than enough room for\/ $G-e$ in\/ $K-v$
 for the extension of\/ $e$ to a rainbow copy of\/ $G$.
\erm

\section{Applications}   \label{s:applic}

Here we apply the results of Sections \ref{s:LB} and
 \ref{s:UB} on various interesting graph classes $\cF$.
In most of the proofs we indicate only the reason why the
 corresponding theorems are applicable.

\btm {\bf\itshape (k-Colorable Decomposition)}   \label{t:kCOL}
 \
Let\/ $\cF_k$ be the class of graphs with chromatic number
 at most\/ $k$ $(k\geq 2)$, and\/ $G$ any graph.
 \begin{itemize}
  \item[$(i)$] If\/ $\chi(G)>2k$, then\/ $f(n, G|\cF_k) =
  (1+o(1)) \, \ex(n,K_{\lceil \chi(G)/k \rceil})$.
  \item[$(ii)$] If\/ $k+1\leq \chi(G)\leq 2k$, then\/
   $f(n, G|\cF_k) = o(n^2)$ and\/ $f\to\infty$ as\/ $n\to\infty$.
 \item[$(iii)$] Otherwise\/ $f(n, G|\cF_k) = 1$.
 \end{itemize}
\etm

\bpf
The chromatic number function $\chi$ is submultiplicative;
 that is, the inequality $\chi(H'\cup H'') \leq \chi(H') \cdot \chi(H'')$
 is valid for any two graphs $H',H''$ on the same vertex set.
As a consequence, $D(K_p |\cF_k)$ contains only graphs with
 chromatic number at least $\lceil \chi(G)/k \rceil$.
On the other hand, the reduced chromatic number
 $\chi_{_{\cF_k}}(G)$ is not larger than $\lceil \chi(G)/k \rceil$.
Indeed, take any optimal proper vertex coloring
 $\vp:V(G)\to\{1,\dots,\chi(G)\}$ of $G$.
Denoting $V_i=\{v\in V(G) : \vp(v)=i\}$, the sets
 $X_t:=\bigcup_{j=1}^k V_{tk+j}$ induce a $k$-colorable graph
 for every $t\geq 0$, while the graph induced by all the edges
 connecting distinct parts $X_{t'},X_{t''}$
 ($t'>t''\geq 0$) has chromatic number
 $\lceil \chi(G)/k \rceil$.
And, of course, if $\chi(G)\leq k$, then $G\in\cF_k$.
\epf

For the next result, recall the notation in Example \ref{exs:F}.

\btm {\bf\itshape (Proper Edge Coloring / Local \Arr\ Number)} 
Let\/ $\cF = \{ tK_2 : t \geq 1 \}$ be the family of matchings.
 \begin{itemize}
   \item[$(i)$] For\/ $p \geq 5$ we have\/
   $\lar(n,K_p) = f(n,K_p |\cF)  = (1+o(1)) \, \ex(n, K_{\lceil p/2 \rceil})$.
   \item[$(ii)$] For every\/ $q \geq 2$ we have\/
    $\lar(n,K_{2q}) - \lar(n,K_{2q-1}) = o(n^2)$.
   \item[$(iii)$] $\lar(n,K_4) = o(n^2)$.
   \item[$(iv)$] For every wheel graph\/ $W_k = C_k + K_1$
    with\/ $k\geq 4$, we have\/
 $\lar(n,W_k) = (1+o(1)) \, n^2\!/4$.
 \end{itemize}
\etm

\bpf
Concerning complete graphs,
 removing a largest matching from $K_p$, we obtain
 $\chi(K_p - \lfloor p/2 \rfloor K_2) = \lceil p/2 \rceil$,
 which is 2 if $p=4$ and at least 3 if $p\geq 5$.
This settles $(i)$--$(iii)$.
Concerning wheels with $k>3$, it suffices to observe that
 removing a matching of $t$ edges can destroy at most $t+1$
 of the $k$ triangles, and the matching number of $W_k$ is
 $\lceil k/2 \rceil < k-1$, hence $\chiF(W_k)=3$.
\epf

For the particular case of $W_3 \cong K_4$, the stronger
 result $\lar(n,K_4)=\Theta(n\sqrt{n})$ is proved in the
 companion paper \cite{ar-M}, using different methods.

In particular, little more than $\ex(n, K_q)$ colors for
 $q\geq 3$ force the presence of a properly colored $K_{2q}$.
As a comparison, little more than $\ex(n, K_q)$ colors force
 the presence of a rainbow $K_{q+1}$ only,
 as proved in \cite{ESS-75}.

In some cases, the bound $o(n^2)$ in Theorem \ref{t:subex} $(ii)$
 can be improved substantially.
Such examples are given in the next result.
A book with $t$ pages is defined as the graph
 $B_t := K_2 + tK_1$.
A nice similar structure is $M_{s,t} := sK_2 + tK_1$
 (hence $B_t=M_{1,t}$).
Even more generally, for any graph $H$ and any positive
 integer $r$, let us denote by $B(H,r)=H+rK_1$ the
 complete join of $H$ and an independent set of size $r$.

\btm {\bf\itshape (Generalized Book Graphs)} 
\
\begin{itemize}
 \item[$(i)$] If\/ $t=1$ or\/ $t=2$, then\/ $\lar(n,B_t)=\Theta(n)$.
 \item[$(ii)$] If\/ $t\geq 3$, then\/ $\lar(n,B_t)=\Theta(n^{3/2})$.
 \item[$(iii)$] For every\/ $s$, if\/ $t\geq t_0(s)$ is sufficiently large, then\/ $\lar(n,M_{s,t})=\Theta(n^{2-1/2s})$.
 \item[$(iv)$] If\/ $|V(H)|=s$, $\Delta(H)\leq t$, and\/ $r\geq r_0(s)$ is sufficiently large,
  then\/ $\lar(n,t,B(H,r))=\Theta(n^{2-1/s})$.
\end{itemize}
\etm

\bpf
We have $B_1=K_3$ with $\lar(n,K_3)=n$ as proved in \cite{ESS-75},
 and $B_2=K_4-e$ with $\lfloor (3n-1)/2 \rfloor \leq \lar(n,K_4-e) \leq 2n-3$ as proved in \cite{ar-M}.
 These results settle $(i)$.

The other three parts are of increasing strength,
 $(ii)$ is the particular case $s=1$ of $(iii)$ provided that
 we prove $t_0(1)=3$, and $(iii)$ follows from $(iv)$
 by putting $r=1$ and $H=sK_2$.
For this reason, the proof of these assertions is very similar.

As it has already been noted, $\lar(n,G)$ corresponds to the
 hereditary family $\cF$ of matchings.
Then the $\cF$-deck of $B_t$ consists of three graphs: $K_{2,t}$,
 the graph obtained from $B_{t-2}$ by attaching two disjoint
 pendant edges to the central edge of the triangles, and
 $B_{t-1}$ with a pendant edge.
Consequently $\lar(n,B_t)\geq \ex(n,\{C_3,C_4\}) + 2$,
 and here the \tur\ function is bounded from below by the
 largest number of edges in a $C_4$-free bipartite graph
 of order $n$, the simplest case of Zarakiewicz's problem.
An asymptotically tight formula for the latter is
 $(\frac{1}{2\sqrt{2}}-o(1)) n^{3/2}$ (see \cite{KST-54}).

In a similar way, if $t>s$, the only $K_3$-free graph in
 the $\cF$-deck of $M_{s,t}$ is $K_{2s,t}$.
Then, since every graph contains a bipartite subgraph with
 at least half of its edges, $\ex(n,D(G|\cF)$ cannot be
 smaller than $\frac{1}{2}\ex(K_{2s,t})$.
Also, for $\lar(n,t,B(H,r))$ and $r$ large enough, the
 hereditary family $\cF$ consists of some 3-chromatic graphs
 and some bipartite ones, where all the latter contain the
 complete bipartite graph $K_{2s,t-rs}$ as a subgraph.
Consequently $\ex(n,D(G|\cF)$ cannot be smaller than half of
 $\ex(K_{2s,t-rs})$ whenever $t\geq (r+2)s$.
Now, the lower bounds for both $(iii)$ and $(iv)$ follow by
 the Koll\'ar--R\'onyai--Szab\'o theorem \cite{KRS-96}, which
 asserts that $\ex(n, K_{a,b}) > c_a\,n^{2-1/a}$ holds
 for every $a>1$, if $b > a!$\,;
in fact also $b > (a - 1)!$ is sufficient for the same
 conclusion, as proved in \cite{ARS-99}.

For upper bounds we quote again \cite{KST-54} where the proof
 of $\ex(n,K_{a,a})=O(n^{2-1/a})$ is well-known to yield
 $\ex(n, K_{a,b}) \leq (1/2 + o(1))\left(b - 1\right)^{1/a}\,n^{2-1/a}$, an
 estimate with the same order of magnitude.
In $(ii)$ and $(iii)$ we take $a=2$ and $a=2s$, respectively,
 with $b=a+s$, and if the number of colors exceeds $\ex(n,K_{a,b})$,
 we find a rainbow $K=K_{a,b}$.
Then in the smaller vertex class of $K$ we select $s$ disjoint
 edges, and omit the at most $s$ vertices from the other class
 as done in the proof of Theorem \ref{l:RCL}.
More generally, for $(iv)$, a suitable choice of the
 parameters is $a=s$, $b=r+|E(H)|\leq r+st/2$.
The condition $\Delta(H)\leq t$ ensures that any coloring
 of a copy of $H$ in the smaller class of $K_{a,b}$ is
 decomposed into monochromatic members of $\cF$.
\epf

The (\emph{linear}) \emph{arboricity} of a graph $G$ is the smallest number
 of (linear) forests $E(G)$ can be decomposed into.

\btm {\bf\itshape (Arboricity / Linear Arboricity)} 
If\/ $\cF$ is the family of forests or the family of
 linear forests, then\/ for every\/ $p \geq 5$ we have\/
   $f(n,K_p |\cF)  = (1+o(1)) \, \ex(n, K_{\lceil p/2 \rceil})$.
\etm

\bpf
Removing any forest (or even any bipartite graph) from $K_p$,
 the chromatic number does not become smaller than $\lceil p/2 \rceil$.
Already a matching verifies that this lower bound is tight.
\epf

For an integer $d\geq 1$ a graph $G$ is called
 \emph{$d$-degenerate} if every subgraph of $G$ contains
 a vertex of degree at most $d$.

\btm {\bf\itshape (d-Degenerate Coloring)} 
For a fixed integer\/ $d\geq 1$ let\/ $\cF$ be the family of\/
 $d$-degenerate graphs. Then\/
 for every\/ $p \geq 2d+3$ we have\/
   $f(n,K_p |\cF)  = (1+o(1)) \, \ex(n, K_{\lceil p/(d+1) \rceil})$.
\etm

\bpf
Since every $d$-degenerate graph is $(d+1)$-colorable,
 $p>2d+2$ implies $\chiF(K_p)\geq \lceil p/(d+1) \rceil >2$.
Tightness is shown by removing
 $ \lfloor\frac{p}{d+1} \rfloor K_{d+1}
   \cup K_{p \, \mathrm{mod} \, (d+1)} $.
\epf

Concerning the third part of Example \ref{exs:F} we have:

\btm {\bf\itshape (t-Defective Decomposition)} 
For a fixed integer\/ $t\geq 1$ let\/ $\cF$ be the family of
 graphs having maximum degree at most\/ $t$.
If $\chi(G) \geq 2t+3$, then\/
   $f(n,G |\cF) = \lar(n,t,G)
     = (1+o(1)) \, \ex(n, K_{\chiF(G)})$.
\etm

\bpf
If $\Delta(D) \leq t$, then $\chi(D) \leq t+1$, hence
 $\chi(G)>2t+2$ implies $\chiF(G)>2$.
\epf

In her recent dissertation, Isabel Harris \cite{H-24t} studied
 the following variant of \arr\ numbers.
Given an integer $k$ and two graphs $G$ and $H$, determine the
 largest possible number of colors in an
 edge coloring of $H$ such that every copy of $G$ in $H$
 necessarily contains more than $k$ edges from a color class.
Restricting attention to $H=K_n$, and to fit better with our
 overall formalism, we let $\arra_k(n,G)$ denote the smallest
 $q$ such that every edge coloring of $K_n$ with at least $q$
 colors contains a copy of $G$
 in which no color occurs more than $k$ times.

\btm {\bf\itshape (k-Anti-Ramsey Numbers)}   \label{t:subram}
Let\/ $\cF_k$ denote the family of graphs that have at most\/
 $k$ edges.
  \begin{itemize}
   \item[$(i)$] If a graph\/ $G$ cannot be made bipartite
    by the omission of at most\/ $k$ edges, then\/
  $\arra_k(n,G)=(1+o(1)) \, \ex(n, K_{\chiFk\!(G)})$.
   \item[$(ii)$] If\/ $\chi(G)=q\geq 3$ and the complete\/
    $q$-partite graph\/ $K_{k+1,\dots,k+1}$ is a subgraph\/
    of\/ $G$, then\/ $\arra_k(n,G) =
      (1+o(1)) \, \ex(n, K_q) =
       (1+o(1)) \, \frac{k-2}{2k-2}\, n^2.$
  \end{itemize}
 
\etm

\bpf
The condition in $(i)$ on $G$ implies $\chiFk\!(G)\geq 3$.
Moreover, the assumption $K_{k+1,\dots,k+1}\subset G$ in $(ii)$
 ensures that $G$ is stable with respect to $\cF_k$.
\epf

In the next theorem we collect results in which a complete
 subgraph of given size has to be found such that all its
 monochromatic edge classes are required to be planar,
 or even more restrictively outerplanar. 

\btm {\bf\itshape (Planar / Outerplanar Decomposition)}   \label{t:PLA}
\
\begin{itemize}
 \item[$(i)$] For\/ $p \geq 9$, $f(n, K_p|\mathit{Planar}) =
  (1+o(1)) \, \ex(n,K_{\lceil p/4 \rceil})$.
 \item[$(ii)$] For\/ $p \geq 7$, $f(n, K_p|\mathit{Outerplanar}) =
  (1+o(1)) \, \ex(n,K_{\lceil p/3 \rceil})$.
 \item[$(iii)$] For\/ $5 \leq p \leq 8$, $f(n, K_p|\mathit{Planar}) =
  o(n^2)$ and\/ $f\to\infty$ as\/ $n\to\infty$.
 \item[$(iv)$] For\/ $4 \leq p \leq 6$, $f(n, K_p|\mathit{Outerplanar}) =
  o(n^2)$ and\/ $f\to\infty$ as\/ $n\to\infty$.
 \item[$(v)$] Otherwise\/ $f(n, K_p|\mathit{Planar}) =
   f(n, K_p|\mathit{Outerplanar}) = 1$.
\end{itemize}
\etm

\bpf
Since $\chi$ is submultiplicative, the Four Color Theorem
 implies that $D(K_p |Planar)$ contains only graphs with
  chromatic number at least $\lceil p/4 \rceil$, and this
 lower bound is tight, as realized by deleting the edges of
 $\lfloor p/4 \rfloor$ vertex-disjoint $K_4$ subgraphs
 (and a further $K_2$ or $K_3$ if $n\equiv 2,3$ (mod 4)).

Similarly but more simply, as every outerplanar graph is
 2-degenerated and hence 3-colorable, we see that
 $D(K_p |Outerplanar)$ contains only graphs with
  chromatic number at least $\lceil p/3 \rceil$, and this
 lower bound is tight, as realized by deleting the edges of
 $\lfloor p/3 \rfloor$ vertex-disjoint triangles
 (and a further edge if $n\equiv 2$ (mod 3)).
\epf

There are cases where $\cF$ is hereditary but $\chiF(G)$ is
 hopeless to determine exactly, yet some meaningful estimates
 can be given.
A prototype of this situation occurs with families defined by
 ``not containing a certain graph $G$ or a family of graphs''
 if chromatic number in $\cF$ is unbounded.   
Then, to obtain some concrete information concerning
   $f(n, G |F) = (1+o(1)) \, \ex(n, K_{\chiF(G)})$,
we may use some estimates how fast the chromatic number
 can grow in $\cF$ with respect to the order of its members.
As a characteristic example, we consider here the class of
 triangle-free graphs, which is hereditary and contains
 graphs with arbitrary large chromatic number.
The determination of the corresponding $\chiF(K_p)$ is in
 close relation with the Ramsey numbers $R(3,p)$, too.

\btm {\bf\itshape (Triangle-Free Decomposition)}   \label{t:K3free}
Let\/ $\cF$ be the class of triangle-free graphs. Then\/
 $$
   f(n, K_p | K_3\mathit{-free}) =
    (1+o(1)) \, \ex(n, K_{\chiF(G)})
     \geq (1+o_n(1)) \, \ex(n, K_{(1/2 +o_p(1)) \sqrt {p \log p} }).
 $$
\etm

\bpf
Davies and Illingworth \cite{DI-22} proved that the
 chromatic number of tri\-angle-free graphs of order $p$ is
  at most $(2+o(1)) \sqrt {p\, / \log p}$ as $p\to\infty$.
Thus, if $D\in\cF$, then
 $\chi(K_p-D) \geq (1/2 +o(1)) \sqrt {p \log p}$.
\epf

We close this section with some consequences of Theorem \ref{t:stable}.
Let us denote by $\chi_t(G)$ the $t$-defective
 chromatic number of $G$, defined as the smallest $k$
  such that there exists a vertex coloring of $G$ where each
 color class induces a subgraph of maximum degree at most $t$
 (admitting also the case $t=0$).
Note that the inequality chain
 $$
   \chi(G) = \chi_0(G) \geq \chi_1(G) \geq \cdots
    \geq \chi_{\Delta(G)}(G) = 1
 $$
  is valid for every graph $G$, directly seen by definition.
The corresponding local \arr\ functions are defined in
 parts 2 and 3 of Example \ref{exs:F}.

A very general implication of Theorem \ref{t:stable} is the following result.

\btm {\bf\itshape (Local \Arr\ Stability)}   \label{t:lar-stab}
Let $t\geq 1$ be any integer.
If\/ $\chi_t(G)=\chi(G)\geq 3$, then\/
 $\lar(n,t,G) = (1+o(1)) \, \ex(n, K_{\chi_t(G)})$.
In particular, if the removal of any matching does not decrease
  the chromatic number of a non-bipartite graph\/ $G$, then\/
 $\lar(n,G) = (1+o(1)) \, \ex(n, K_{\chi(G)})$.
\etm

\bpf
It suffices to take $\cF$ as the family of graphs with
 maximum degree at most $t$, and for $t=1$ the family of matchings.
\epf

Finally, we provide some existence results.

\btm {\bf\itshape (Planar / Outerplanar Quadratic Growth)}   \label{t:quad}
\
Let\/ $\cF = \{ tK_2 : t \geq 1 \}$ be the family of matchings.

\begin{itemize}
 \item[$(i)$] There exist planar graphs\/ $G$ with\/
  $\lar(n,G) > n^2\!/3$.
 \item[$(ii)$] There exist outerplanar graphs\/ $G$ with\/
  $\lar(n,G) > n^2\!/4$.
\end{itemize}
\etm

\bpf
Cowen, Cowen, and Woodall proved in \cite{CCW-86} that
 there exist planar graphs $G$ not admitting
 any 1-defective vertex coloring with three colors;
 and also outerplanar graphs $G'$ that do not admit
 any 1-defective vertex coloring with two colors.
All such graphs have $\chi_1(G) = \chi(G) = 4$ and
 $\chi_1(G') = \chi(G') = 3$, respectively.
As a consequence, those $G$ and $G'$ are stable
 with respect to $\cF$, implying
 $f(n,G|\cF) > \ex(n,K_4) = \lfloor n^2\!/3 \rfloor$
  and $f(n,G'|\cF) > \ex(n,K_3) = \lfloor n^2\!/4 \rfloor$.
\epf

\btm {\bf\itshape (High-Girth Quadratic Growth)}   \label{t:girth}
\
For any fixed integer\/ $t\geq 1$ let\/ $\cF$ be the family of
 graphs having maximum degree at most\/ $t$.
For any\/ $k\geq 2$ and any\/ $g\geq 3$
 there exist\/ $k$-degenerate graphs\/ $G$ such that\/
 the girth of\/ $G$ is at least\/ $g$, and\/
  $f(n,G|\cF) = (1+o(1)) \, \ex(n,K_{k+1})
   \approx \frac{k-1}{2k} n^2$.
\etm

\bpf
Kostochka and Ne\v set\v ril
 \cite{KN-99} proved, for any $k$, the existence of
 $k$-degenerate graphs $G$ with girth at least $g$
 that are not $k$-colorable with defect $t$.
This means $\chi_t(G) = \chi(G) = k+1$, due to the fact
 $\chi_t(G) \leq \chi(G) \leq k+1$.
Thus, $G$ is stable with respect to $\cF$.
\epf

We note that there is a wide range of further natural families $\cF$ admitting
 infinitely many graphs that are stable with respect to $\cF$. This issue will be discussed in detail in the companion paper \cite{ar-M}.

\section{Concluding remarks and open problems}
\label{s:concl}

In this paper we introduced a very general extremal problem
 concerning edge decompositions of graphs.
Many extensively studied areas turn out to be just particular
 cases of our model.
As is well known, \tur-type extremal graph theory makes a
 substantial distinction between bipartite and non-bipartite graphs.
We have proved that the same phenomenon applies also in our model.
In particular, in a well-defined way of non-bipartiteness
 with respect to any given hereditary family of graphs,
 we have derived a tight formula of the asymptotic growth
 of the decomposition function $f(n,G|\cF)$ under study.

Exact results for particular types of graphs $G$
 and graph classes $\cF$ would be of definite interest.
Moreover, the
``degenerate'' (bipartite) cases still offer an open area
 to explore, although the corresponding upper bound of type
 $o(n^2)$ has also been proved here, and even stronger
 asymptotics have been determined in some cases.

In summary:
Much like in classical \tur\ and \arr\ problems, the main
 future challenges in this area can be considered in two
  directions:

\bpm
Let\/ $\cF$ be a hereditary family of graphs, and\/ $G$ a
 given graph.
 \begin{itemize} 
  \item[$(i)$] Find better asymptotic with a smaller concrete
   error term for\/ $f(n,G|F)$ in case\/ $\chiF(G) \geq 3$.
  \item[$(ii)$] Find better lower and upper bounds for\/
   $f(n,G|F)$ in case\/ $\chiF(G) = 2$.
 \end{itemize}
\epm

\subsection{Tightness of the Independent Transversal Lemma}
\label{ss:c-ITL}

In connection with Theorem \ref{l:ITL} $(i)$, let us introduce
 the notation $s(m,d)$ for the smallest integer $s$
 with the following property:
If $D$ is a balanced multipartite digraph with $m$ vertex
 classes $V_1,\dots,V_m$, each of size $s$, and $D$ has
 maximum out-degree $\Delta^+(D)\leq d$, then there exists an
 independent transversal in $D$, that is an independent set
  meeting every class $V_i$.
The general upper bound
 $$
   s(m,d) \leq md + 1
 $$
  is valid by Theorem \ref{l:ITL}.
Although this estimate holds with strict inequality in
 most of the cases, it is not far from being tight, and
 for infinitely many values of $m$ and $d$ it is best possible.

In this final subsection we characterize the case of equality,
 and also prove a general lower bound on $s(m,d)$.
We postponed this analysis here because the actual value of
 $s(m,d)$ is irrelevant in connection with applicability
 in the main subject, $f(n,G|\cF)$, of this paper.

\btm
The following estimates are valid on\/ $s(m,d)$.
  \begin{itemize}
   \item[{$(i)$}] For every\/ $m\geq 3$ and\/ $d\geq 1$,
    $s(m,d) \geq md - d + 1$; and if\/ $m\leq d+1$, then\/
     $s(m,d) \geq md - d + m$.
   \item[{$(ii)$}] The equality\/ $s(m,d) = md + 1$ holds
     if and only if\/ $d$ is divisible by\/ $m-1$
     (including\/ $m=2$ for all\/ $d$).
   \item[{$(iii)$}] In particular, $s(m,1)=m$ for every\/
    $m\geq 3$, moreover\/ $s(2,d)=2d+1$ for every\/ $d\geq 1$,
     and\/ $s(d+1,d)=d^2+d+1$ for every\/ $d\geq 2$.
  \end{itemize}
\etm

\bpf
$(i)$\quad 
The following construction shows $s(m,d) > (m-1)d$.
Partition each $V_i$ into $m-1$ subsets,
 $V_i=A_{i,1}\cup\cdots\cup A_{i,m-1}$ with $|A_{i,j}|=d$
 for $i=1,\dots,m$ and $j=1,\dots,m-1$.
    For $i,j=1,\dots,m-1$, orient $d$ arcs
     from each $v\in A_{i,j}$ to $A_{m,i}$.
Then any selection $v_1,\dots,v_{m-1}$ of vertices $v_i\in V_i$
 has its out-neighborhood being the entire $V_m$, hence
 an independent transversal does not exist.

If $m\leq d$, then we construct a digraph with $m$ vertex
 classes $V_1,\dots,V_m$ with $s=(m-1)(d+1)$ vertices in each part.
Let $V_m = V_{m,1} \cup \cdots \cup V_{m,m-1}$ with
 $|V_{m,i}|=d+1$ for all $1\leq i\leq m-1$,
  $V_{m,i} = \{v_{m,i,1},\dots,v_{m,i,d+1}\}$.
For each $i\leq m-1$, partition $V_i$ into $d+1$ subsets,
 $V_i=A_{i,1}\cup\cdots\cup A_{i,d+1}$, with $|A_{i,j}|=m-1$.
For $i=1,\dots,m-1$ and $j=1,\dots,d+1$, orient $d$ arcs
 from each $v\in A_{i,j}$ to $V_{m,i}\smin\{v_{m,i,j}\}$,
 and $m-1$ arcs from $v_{m,i,j}$ to $A_{i,j}$.
Hence all out-degrees in $V_m$ are $m-1\leq d$, and the
 out-degree of every vertex not belonging to $V_m$ is $d$.

Let $T$ be any set meeting each $V_i$ in one vertex.
Say, $T\cap V_m$ is the vertex $v=v_{m,1,1}$.
Consider now the vertex $w$ of $T\cap V_1$.
If $w\in A_{1,1}$, then $\overrightarrow{vw}$ is an arc.
If $w\in A_{1,j}$ for some $j>1$, then $\overrightarrow{wv}$ is an arc.
Thus, $T$ cannot be an independent transversal.

\msk

$(ii)$\quad 
First we provide constructions showing $s(m,d) > md$ in the
 given cases.
Denote $k:=\frac{d}{m-1}$, which is now supposed to be
 an integer.
We need to construct $m$-partite digraphs with out-degree
 $(m-1)k$, and with $m(m-1)k$ vertices in each part $V_i$
 ($i=1,\dots,m$), not admitting any independent transversals.
We define a vertex partition as follows.
 \begin{itemize}
  \item Let $V_i=V_{i,1}\cup\cdots\cup V_{i,m}$,
   with $|V_{i,j}| = (m-1)k = d$ for all
   $1\leq i\leq m-1$ and all $1\leq j\leq m$.
  \item Let $V_m=V_{m,1}\cup\cdots\cup V_{m,m-1}$,
   with $|V_{m,j}| = mk = d+k$ for all
   $1\leq j\leq m-1$; and refine further each $V_{m,j}$ as
   $V_{m,j}=V_{m,j,1}\cup\cdots\cup V_{m,j,m}$ with
   $|V_{m,j,l}| = k$ for all $1\leq j\leq m-1$ and all $1\leq l\leq m$.
 \end{itemize}
Now the arc set of the digraph is defined as
 \begin{itemize}
  \item $d$ arcs from each $v\in V_{m,j,l}$ to the entire
   $V_{j,l}$,
  \item $d$ arcs from each $w\in V_{j,l}$ to the entire
   $V_{m,j}\smin V_{m,j,l}$.
 \end{itemize}
Then the underlying undirected graph induced by any
 $V_i\cup V_{m,i}$ ($1\leq i\leq m-1$) is a complete
  bipartite graph $K_{m(m-1)k,mk}$.
Hence the neighborhood of any vertex in $V_m$ is an entire
 vertex class, preventing the digraph from admitting an
 independent transversal.

A transparent particular structure is obtained for $m=2$,
 where the notation can be simplified to $V_i=A_i\cup B_i$
  with $|A_i|=|B_i|=d$ for $i=1,2$, and orienting
 $d$ arcs from each $v\in A_1$ to $A_2$,
 from each $v\in B_1$ to $B_2$, from each 
 $v\in A_2$ to $B_1$, and from each $v\in B_2$ to $A_1$.
The underlying undirected graph of those arcs is $K_{2d,2d}$.

\ssk

It remains to prove that in case $(m-1) \nmid d$, if
 $D$ is a digraph in which all the
 $m$ classes $V_i$ are of size $md$, and all out-degrees are
 at most $d$, then there exists an independent transversal.
Following the proof of Theorem \ref{l:ITL} $(i)$, each of the
 $m^2d^2$ arcs alone destroys exactly $(md)^{m-2}$ of the
 $(md)^m$ candidate sets for independent transversals.
Hence, to destroy all, it would be necessary that none of
 those $m$-element sets is counted more than once.
This condition excludes the following substructures:
 \begin{enumerate}
  \item a directed 2-cycle $C_2$ in $D$,
   that is two oppositely oriented arcs $\overrightarrow{uv}$
   and $\overrightarrow{vu}$;
  \item two disjoint arcs (an oriented $2K_2$) that together
   meet four vertex classes;
  \item two adjacent arcs (an oriented $P_3$) that
    together meet three vertex classes.
 \end{enumerate}
The proof of $(ii)$ will be done if we show that at least one
 of these three ``forbidden'' situations must occur
 whenever $m\geq 3$ and $(m-1) \nmid d$.

Construct an auxiliary undirected graph $H$ of order $m$, with
 vertices $x_1,\dots,x_m$.
Draw an edge $x_ix_j$ if there is at least one arc between
 $V_i$ and $V_j$ in $D$.
The exclusion of $2K_2$ means that H is either a star
 $K_{1,m-1}$ or a triangle $K_3$.
The latter is possible only if $m=3$, because all out-degrees
 in $D$ are positive.

Assume first that $H$ is a star, say with center $x_m$.
Then $V_m$ can be partitioned as $A_1\cup\cdots\cup A_{m-1}$
 such that all arcs between $V_i$ and $V_m$ meet $A_i$,
 for every $i<m$.
Suppose there is a non-edge $v_1v_m$, for some $v_1\in V_1$
 and $v_m\in A_m$.
Then this vertex pair can be
 supplemented to an independent transversal by selecting an
 arbitrary vertex $v_i$ from each $V_i$, $1<i<m$.
Hence we can assume that each $V_i\cup A_i$ induces a
 complete bipartite graph $G_i\cong K_{md,q_i}$ in the
 undirected underlying graph of $D$.
We also have $q_1+\ldots+q_{m-1}=|V_m|=md$.
To have a complete bipartite graph it is necessary that all
 possible arcs are present in at least one direction, and
 in fact in exactly one direction, as $C_2$ has been excluded:
  $$
    ( md + q_i ) d = d|A_i| + d|A_m| = |A_i|\cdot|A_m|
      = mdq_i \,,
  $$
  $$
    q_i = \frac{md}{m-1} \,.
  $$
Since $\gcd(m-1,m)=1$, it follows that $m-1$ is a divisor of $d$.

Assume now that $m=3$ and $H$ is a triangle.
The exclusion of $P_3$ yields a partition $A_i\cup B_i=V_i$
 for $i=1,2,3$ such that each arc is included in one of the
 sets $A_1\cup B_2$, $A_2\cup B_3$, $A_3\cup B_1$,
 none of them being edgeless.
Choosing any vertices $v_1\in A_1$, $v_2\in A_2$, $v_3\in A_3$,
 an independent transversal is obtained.

\msk

\nin
$(iii)$\quad 
The cases of $s(2,d)$ and $s(d+1,d)$ satisfy the arithmetic
 conditions sufficient for $s(m,d)=md+1$ described in $(ii)$.
On the other hand, divisibility does not hold for $s(m,1)$
 if $m>2$, hence $s(m,1)\leq m$ by $(ii)$, while we also have
 $s(m,1)\geq m$ by $(i)$.
\epf

\bpm
Complete the determination of\/ $s(m,d)$, or a least improve the estimates.
\epm


\begin{thebibliography}{99}

\bibitem{ADK-24+}
N. Alon, C. Defant, N. Kravitz:
Rainbow stackings of random edge-colorings.
arxiv:2405.14795 

\bibitem{ARS-99}
N. Alon, L. R\'onyai, T. Szab\'o:
Norm-graphs: variations and applications.
J. Combinatorial Theory, Ser. B 76 (1999), 280--290.

\bibitem{A-etal-24}
P. Ara\'ujo, T. Martins, L. Mattos, W. Mendon\c ca,
 L. Moreira, G. O. Mota:
On the anti-Ramsey threshold for non-balanced graphs.
 Electronic J. Combinatorics 31:1 (2024), \#P1.70.

\bibitem{AC-24}
M. Axenovich, F. C. Clemen:
Rainbow subgraphs in edge-colored complete graphs: 
   Answering two questions by Erd\H os and Tuza.
J. Graph Theory 106:1 (2024), 57--66.

\bibitem{B-etal-24+}
N. Behague, R. Hancock, J. Hyde, S. Letzter, N. Morrison:
Thresholds for constrained Ramsey and anti-Ramsey problems.
arXiv:2401.06881

\bibitem{BGR-15}
A. Bialostocki, S. Gilboa, Y. Roditty:
Anti-Ramsey numbers of small graphs.
Ars Combinatoria 123 (2015), 41--53.

\bibitem{hered-surv}
M. Borowiecki, I. Broere, M. Frick, P. Mih\'ok, G. Semani\v sin:
A survey of hereditary properties of graphs.
Discussiones Mathematicae Graph Theory 17:1 (1997), 5--50.

\bibitem{CPS-22}
Y. Caro, M. Petru\v sevski, R. \v Skrekovski:
Remarks on odd colorings of graphs.
Discrete Applied Mathematics 321 (2022), 392--401.

\bibitem{CPS-23}
Y. Caro, M. Petru\v sevski, R. \v Skrekovski:
Remarks on proper conflict-free colorings of graphs.
Discrete Mathematics 346 (2023), 113221.

\bibitem{ar-M}
Y. Caro, Zs. Tuza:
Coloring problems with parity constraints, inspired
  by anti-Ramsey theory and the odd-coloring problem.
Manuscript, 2024, in preparation.

\bibitem{CCW-86}
L. J. Cowen, R. H. Cowen, D. R. Woodall:
Defective colorings of graphs in surfaces:
Partitions into subgraphs of bounded valency.
Journal of Graph Theory 10:2 (1986), 187--195.

\bibitem{DI-22}
E. Davies, F. Illingworth:
The $\chi$-Ramsey problem for triangle-free graphs.
SIAM J. Discret. Math. 36:2 (2022), 1124--1134.

\bibitem{ES-66}
P. Erd\H os, M. Simonovits:
 A limit theorem in graph theory.
Studia Sci. Math. Hungar. 1 (1966), 51--57.

\bibitem{ESS-75}
P. Erd\H os, M. Simonovits, V. T. S\'os:
 Anti-Ramsey theorems.
In: Infinite and Finite Sets,
Proc. colloq. dedicated to P. Erd\H os on his 60th birthday,
 Keszthely, 1973;
  Colloq. Math. Soc. J\'anos Bolyai, 10, Vol. II, 633--643, North-Holland, Amsterdam, 1975.

\bibitem{ES-46}
P. Erd\H os, A. H. Stone:
 On the structure of linear graphs,
Bull. Amer. Math. Soc. 52 (1946), 1087--1091.

\bibitem{dyn-surv}
S. Fujita, C. Magnant, Y. Mao, K. Ozeki:
Rainbow generalizations of Ramsey theory -- A dynamic survey.
Theory and Applications of Graphs Vol. 0: Iss. 1 (2014)
DOI: 10.20429/tag.2014.000101 
https://digitalcommons.georgiasouthern.edu/tag/vol0/iss1/1

\bibitem{GS-22}
S. Glock, B. Sudakov:
An average degree condition for independent transversals.
J. Combinatorial Theory, Ser. B 154 (2022), 370--391.

\bibitem{HW-23+}
P. Haxell, R. Wdowinski:
Degree criteria and stability for independent transversals.
arXiv:2305.10595

\bibitem{H-24t}
I. Harris:
Avoiding $k$-Rainbow Graphs in Edge Colorings of $K_n$
 and other Families of Graphs.
PhD Thesis, Auburn university, AL, 2024.
https://etd.auburn.edu/bitstream/handle/10415/9200/I\%20Harris\%20 Dissertation\%20Final.pdf?sequence=2
\,(accessed on May 8, 2024)

\bibitem{JLZ-19}
Y. Jia, M. Lu, Y. Zhang:
Anti-Ramsey problems in complete bipartite graphs for $t$
 edge-disjoint rainbow spanning subgraphs:
  Cycles and matchings.
Graphs and Combinatorics 35:1 (2019), 1011--1021.

\bibitem{JG-24}
Z. Jin and J. Gu:
Rainbow disjoint union of clique and matching in edge-colored complete graph.
Discussiones Mathematicae Graph Theory 44 (2024), 953--970.

\bibitem{KRS-96}
J. Koll\'ar, L. R\'onyai, T. Szab\'o:
Norm graphs and bipartite Tur\'an numbers.
Combinatorica 16:3 (1996), 399--406.

\bibitem{KN-99}
A. V. Kostochka, J. Ne\v set\v ril:
Properties of Descartes’ construction of triangle-free graphs with high chromatic number.
Combin. Probab. Comput. 8:5 (1999), 467--472.

\bibitem{KST-54}
T. K\H ov\'ari, V. T. S\'os, T. Tur\'an:
On a problem of K. Zarankiewicz.
Colloq. Math. 3 (1954), 50--57.

\bibitem{LSS-21}
Y. Lan, Y. Shi, Z. Song:
Planar Tur\'an number and planar anti-Ramsey number of graphs.
Operations Research Transactions 25:3 (2021), 200--216.

\bibitem{LTY-24+}
T. Li, Y. Tang, G. Yan:
Anti-Ramsey numbers of expansions of doubly edge-critical graphs in uniform hypergraphs.
arXiv:2405.11207

\bibitem{LTJ-09}
X Li, J. Tu, Z. Jin:
Bipartite rainbow numbers of matchings.
Discrete Mathematics 309:8 (2009), 2575--2578.

\bibitem{LS-23+}
X. Liu, J. Song:
Hypergraph anti-Ramsey theorems.
arxiv:2310.01186

\bibitem{WW-22}
I. M. Wanless, D. R. Wood:
A general framework for hypergraph coloring.
SIAM J. Discrete Math. 36:3 (2022), 1663--1677.

\bibitem{W-intro}
D. B. West:
Introduction to Graph Theory.
Prentice Hall, 1996.

\bibitem{XLYX-22}
C. Xiang, Y. Lan, Q. Yan, C. Xu:
The outer-planar anti-Ramsey number of matchings.
Symmetry 14:6 (2022), 1252.

\end{thebibliography}
\end{document}